\documentclass[a4paper,10pt]{article}
\usepackage{amsmath}
\usepackage{amssymb,amsfonts}
\usepackage[english]{babel}
\usepackage{gastex}
\usepackage{graphicx}
\usepackage{cite}
\usepackage[margin=3cm]{geometry}
\usepackage{authblk}

\newtheorem{prop}{Proposition}
\newtheorem{lemm}{Lemma}
\newtheorem{theorem}{Theorem}
\newtheorem{corollary}{Corollary}
\newcommand{\qed}{\hfill\rule{1ex}{1ex}}

\begin{document}

\title{Polynomial identities of the adjoint Lie algebra of $M_{1,1}$}

%Departamento de Matem\'atica,\\
%Universidade de Bras\'\i lia,\\
%70910-900 Bras\'\i lia, DF, Brazil }

\author[1]{Olga Finogenova
\thanks{ Partially supported by Russian Foundation for Basic Research (grants 17-01-
00551, 16-01-00795), by the Ministry of Education and Science of the Russian Federation (project
1.6018.2017/8.9) and the
Competitiveness Enhancement Program of Ural Federal University; 
e-mail \texttt{ob.finogenova@urfu.ru}}
}

\author[2]{ Irina Sviridova 
\thanks{Supported by CNPq, CAPES, FAPESP; e-mail \texttt{I.Sviridova@mat.unb.br}} }

\affil[1]{\small Institute of Natural Sciences and Mathematics, Ural Federal University, Ekaterinburg, Russia\\ 
Department of Mathematics, Technion, Haifa, Israel}

\affil[2]{\small Departamento de Matem\'atica,
Universidade de Bras\'\i lia,
70910-900 Bras\'\i lia, DF, Brazil}

\date{}
\maketitle

%\vspace{1cm}



\vspace{1cm}

\noindent{\bf Abstract}.  We search an identity basis for the adjoint Lie algebra of the algebra $M_{1,1}(K)$ 
over a field, where $K$ is either the infinite generated Grassmann algebra $E$  or $E^1$, the variant of the algebra with $1$.  
In particular, we prove that in the case of an infinite base field of characteristic different from two 
the identities of $M_{1,1}(E^1)$  are exactly all the consequences 
of the identity $[x,y,[z,t],v]=0$. We also find an identity basis of $M_{1,1}(E)$ consisting of three identities.

\vspace{1cm}

\subsection*{Introduction}

Let $\mathbb F$ be a field and let $L(X)$ be the free linear $\mathbb F$-algebra  freely generated by the countable set 
$X = \{ x_1, x_2, \ldots \}$. 
Suppose $A$ is a linear $\mathbb F$-algebra. An element $f(x_1, x_2, \ldots)$ from $L(X)$ is called a {\em polynomial identity} 
of $A$ if $f(a_1, a_2, \ldots)=0$ for all $a_1, a_2, \ldots$ from $A$. If one considers associative algebras, then 
$L(X)$ is the free associative algebra; if it is considered Lie algebras, then  $L(X)$ is the free Lie algebra. 
%An algebra satisfying a non-trivial identity  

  The Grassmann algebra $E$ is the associative $\mathbb F$-algebra generated by a countable set of generators
 $e_1, e_2, \ldots $ with relations $e_ie_j = -e_je_i$ and $e_i^2=0$ for all $i, j$. 
Let  $E^1$  be the Grassmann algebra with 
identity element $1$. We denote by $E_0$ the span of all elements from $E$ of even length, by  $E_1$ the span of 
elements of odd length, and by $E_0^1$ the variant of $E_0$ with $1$. Clearly, $E_0$ consists 
of central elements of $E$ and  for any $x, y$ from $E_1$ we have $xy=-yx$ and $x^2=0$. 
The algebra of all matrices of the kind $\left(
\begin{array}{cc}
a & b \\
d & c
\end{array}\right)$, where $a, c \in E_0$ and $b, d \in E_1$, is denoted by $M_{1,1}(E)$. The algebra 
$M_{1,1}(E^1)$ consists of the same matrices with $a, c \in E_0^1$. 

%The algebra $M_{1,1}(E^1)$ plays an important role in the theory of $PI$-algberas. 
Each of the algebras generates  a ``small'' so-called
{\em verbally prime} varieties. Such varieties play a key role in the theory 
developed by A.~Kemer for the solution of the long standing Specht problem (see, for example,~\cite{Kemer91}). 
Therefore, the identities of $M_{1,1}(E^1)$ and $M_{1,1}(E)$ were investigated 
by many authors from different points of view. A basis of identities of the algebras in the case of a 
zero characteristic base field are found by A.~Popov~\cite{Popov}. The graded identities of $M_{1,1}(E^1)$  are 
investigated by O.~Di Vincenzo~\cite{Vincenzo}. The subvarieties of the variety generated by our algebra are studied 
by L.~Samoilov~\cite{Samoilov}. Some other questions connected with verbally prime varieties
and, in particular, with the algebras, are considered by P.~Koshlukov with
co-authors~\cite{AK06, AFK04, AFK05, KM} and by other researchers. 

Our purpose is to study Lie identities of $M_{1,1}(E)$  and $M_{1,1}(E^1)$. Let us recall that 
with any associative algebra $\langle A,+,\cdot\rangle$, the Lie algebra is associated in a natural way.
The algebra is a so-called {\em adjoint} Lie algebra $\langle A, +, [ , ] \rangle$, where the
multiplication $[ , ]$ is defined by letting
$ [a, b] = a\cdot b - b \cdot a$ for all $a, b  \in A$. We consider the adjoint Lie algebras
both of $M_{1,1}(E)$  and of $M_{1,1}(E^1)$ using the same notations for the Lie algebras as for 
corresponding associative ones. The context will make it clear which algebra, associative or  adjoint Lie, is meant.

    %\vspace{1cm}
We write $[x_1,\ldots, x_n]$ for left-normalized product $[[\ldots [x_1,x_2], \ldots ], x_{n-1}], x_n]$ and 
$[x, z^{(m)}]$
for  $[x, \underbrace{z, \ldots , z}_m]$.

\smallskip

Now we are ready to formulate the main results of the paper.

%\begin{corollary}
\begin{theorem}
\label{tinf}
If $\;\mathbb F$ is an infinite field of characteristic different from 2, then the polynomial identities 
of the Lie algebra $M_{1,1}(E^1)$ admit a basis consisting of the identity $[x,y,[z,t],u]=0$.
\end{theorem}
%\end{corollary}

%\begin{corollary}
\begin{theorem}
\label{t0}
If $\;\mathbb F$ is a zero characteristic field, then  the polynomial identities 
of the Lie algebra $M_{1,1}(E)$ admit a basis consisting of the identity $[x,y,[z,t],u]=0$.
\end{theorem}
%\end{corollary}

\begin{theorem}
\label{Ep}
Let $\mathbb F$ be a field, maybe finite, of characteristic $p > 2$. 
Then the polynomial identities 
of the Lie algebra $M_{1,1}(E)$ admit a basis consisting of the identities 
$[x,y,[z,t],u]=0$, $[x,y, \underbrace{z, \ldots , z}_p, t] =0$, 
$[x,\underbrace{y, \ldots , y}_p,\underbrace{x, \ldots , x}_{p-1},  z] = 0$.
\end{theorem}

The Theorems~\ref{tinf} shows that in the case of infinite field 
the algebra $M_{1,1}(E)$ generates the variety of all centre-by-metabelian 
Lie algebras. Another  algebra as a ``carrier'' of the variety is offered by A.~Krasilnikov~\cite{Kras2008}. 
The questions of having  finite identities bases for some Lie algebras and Lie rings are discussed there.

\subsection*{Proofs of Theorem~\ref{tinf} and Theorem~\ref{t0}}

Let us recall some base notations. For a monomial $h$ and a variable $x$ denote by $deg_x(h)$ the degree 
$h$ in the variable $x$, the number of
times that $x$ occurs in $h$. If all monomials of  a polynomial $f$ have the same degree $m$ in $x$, 
then $f$ is said to be 
{\em homogeneous in}~$x$ of degree $m$ in $x$. In this case we write $deg_x(f) = m$.
A {\em multihomogeneous} polynomial is the polynomial that is homogeneous with respect to all its variables.
The length of monomials of such polinomial $f$ is said to be the {\em degree} of $f$. 

The sum of all monomials of a fix degree in $x$ from $f$ is called a {\em homogeneous in $x$ component} of $f$. If 
one fixes degrees in all variables, then the sum of all monomials with such degrees is called 
a {\em multihomogeneous component} of $f$.  
  
We write $f(\bar x)$ instead of $f(x_1, x_2, \ldots)$ for brevity.

\medskip

\noindent {\bf Remark.} \label{homog} {\em Let $f$ be an identity of an algebra. The following fact is well-known. 
If the maximal degree in $x$ of all monomials from $f$ is less than the base field 
order or the base field is infinite, then the identity $f$ is equivalent to the set of identities, obtained as
homogeneous in $x$ components of $f$.}
   
\medskip

The purpose of the part is to prove the next
two statements.
  
%\begin{theorem}
\begin{prop}
\label{ml}
Let $\mathbb F$ be a field of characteristic different from two. Then every multilinear identity of the 
Lie algebra $M_{1,1}(E)$ or of the Lie algebra $M_{1,1}(E^1)$
is  a consequence of the identity $[x,y,[z,t],u]=0$.  
\end{prop}
%\end{theorem}
 
%\begin{theorem}
\begin{prop}
\label{mh1}
Let $\mathbb F$ be a field of characteristic different from two. Then every multihomogeneous identity 
of the Lie algebra $M_{1,1}(E^{1})$  is 
a consequence of the identity $[x,y,[z,t],u]=0$.  
\end{prop}
%\end{theorem}

The Propositions imply Theorems~\ref{tinf} and~\ref{t0}. Indeed,   
%it is well-known that the identities of algebras over an infinite 
%field are generated by multihomogeneous identities. Thus, 
Theorem~\ref{tinf} is an evident consequence of 
Proposition~\ref{mh1} because of the remark above. 

In the case of zero characteristic field identities of each algebra follow from multilinear ones. 
Hence, to show the truth of Theorem~\ref{t0} is sufficient to notice that the algebra 
$M_{1,1}(E^1)$ satisfies $[x,y,[z,t],u]=0$ and 
contains $M_{1,1}(E)$ as a subalgebra. This fact and Proposition~\ref{ml}  mean that 
$M_{1,1}(E^1)$ and $M_{1,1}(E)$ have  
the same multilinear identities. Theorem~\ref{t0} is proved. \qed

\medskip

There are several methods to prove Propositions~\ref{ml} and~\ref{mh1}. The first one is to  apply
 Kuzmin's result~\cite[Theorem~4]{Kuzmin}. The theorem 
describes the bases of the free centre-by-metabelian Lie algebras. 
The second way is to use the result of I.~Volichenko~\cite[Theorem 1]{Voliche}. He proved that the variety 
of all centre-by-metabelian Lie algebras is a minimal  variety satisfying no identity of  some special kind.   
Also one can use the full description of the lattice of all centre-by-metabelian Lie algebras to prove 
Proposition~\ref{ml}. The description is found by S.~Mishchenko~\cite{Mishch}. 
We prove Propositions by the way different from the ways above. 

\medskip

\iffalse
\begin{corollary}
If $\;\mathbb F$ is a field of characteristic $p>2$, then the polynomial identities 
of the Lie algebra $M_{1,1}(E)$ admit a basis consisting of the identities 
$[x,y,[z,t],u]=0$, $[x,y, \underbrace{z, \ldots , z}_p, t] =0$, 
$[x,\underbrace{y, \ldots , y}_p,\underbrace{x, \ldots , x}_{p-1},  z] = 0$.
\end{corollary}
\fi

From now on, $\mathbb F$  will be a field of characteristic different from two.
Let us denote

\begin{equation}
\label{cma}
cm = [x,y,[u,v],z]
\end{equation}

%${\mathcal{CM}} = \{ [x,y,[z,t],u]\} ^L$

\begin{lemm}
\label{Ins}
%Let $\mathbb F$ be a ring,  $\frac{1}{2} \in F$, and 
Let 
%$char \mathbb F \ne 2$ and let 
$f(\bar x) = \sum_{i j k} \alpha_{i j k} [x_i,x_j,\ldots , x_k]$ 
be a multilinear Lie polynomial with coefficients from $\mathbb F$. Suppose $g(\bar x, y,z) = \sum_{i j k} \alpha_{i j k} [x_i,x_j, y, z, \ldots , x_k] $. Then 
$g$ is a consequence of $f$ and~$cm$.
\end{lemm}

\noindent {\bf Proof}. It is easy to see that modulo~(\ref{cma}) the following holds: 
\begin{align*}
f(x_1, \ldots, & [x_i, y], \ldots) \equiv  \\
& \sum_{j , k \ne i} (\alpha_{i j k} [[x_i,y], x_j \ldots , x_k] + \alpha_{k i j} [x_k, [x_i,y], \ldots , x_j] + \alpha_{j k i} [x_j, x_k, \ldots , [x_i,y]]),
\end{align*}
\begin{align*}
[x_i,x_j,y,z, \ldots , x_k] & = [x_i,[x_j, y, z] \ldots , x_k] + [[x_i,y,z], x_j,\ldots , x_k],  \\
[x_i,x_j, \ldots , y, z, x_k] & = [x_i, x_j,  \ldots , [x_k, y, z]] + [x_i, x_j,\ldots , x_k, z , y]. 
\end{align*}

Hence, we have modulo~(\ref{cma}) 
\begin{align*}
 & 2 g(\bar x, y, z) \equiv  \sum_{i j k} \alpha_{ijk} ([x_i,[x_j, y, z] \ldots , x_k] + [[x_i,y,z], x_j,\ldots , x_k]) +\\ 
& \sum_{i j k} \alpha_{i j k} ([x_i, x_j,  \ldots , [x_k, y, z]] + [x_i, x_j,\ldots , x_k, z , y]) = \\
 & \sum_{i j k} \alpha_{i j k} [x_i,[x_j, y, z] \ldots , x_k]  + \sum_{i j k} \alpha_{i j k} [[x_i,y,z], x_j,\ldots , x_k] + \\
& \sum_{i j k} \alpha_{i j k} [x_i, x_j,  \ldots , [x_k, y, z]] +   \sum_{i j k} \alpha_{ijk}  [x_i, x_j,\ldots , x_k, z , y] = \\
& \sum_{i j k} \alpha_{i j k} [[x_i, y, z], x_j \ldots , x_k] +  \sum_{i j k} \alpha_{k i j} [x_k,[x_i, y, z] \ldots , x_j] +  \\
& \sum_{i j k} \alpha_{j k i} [x_j, x_k,  \ldots , [x_i, y, z]]  + [f(x_1, \ldots), z, y] =\\
& \sum_{i} (\sum_{j , k \ne i} (\alpha_{i j k} [[x_i,y,z], x_j \ldots , x_k] + \alpha_{k i j} [x_k, [x_i,y,z], \ldots , x_j] + \alpha_{j k i} [x_j, x_k, \ldots , [x_i,y,z]])) + \\
& + [f(x_1, \ldots), z, y]  \equiv  \sum_{i} f(x_1, \ldots, [x_i,y,z], \ldots) + [f(x_1, \ldots), z, y].  \;\;\; 
\end{align*}

\noindent {\bf Remark.} {\em Lemma~\ref{Ins} is true for the polynomials with the coefficients from a ring 
contained $\frac{1}{2}$.}  

%$L = \{ [x,y,[z,u],v]\}^L$, 

%$\mathcal{CM} $  is the relatively free algebra of the variety $\var \{ [x,y,[z,u],v] \}$ of Lie F-algebras . 

\begin{corollary}[Folklore]
%well-known (see Popov) %Let $\mathbb F$ be a ring,  $\frac{1}{2} \in F$. Then 
For each integer $k \ge 0$ and variables 
$t_1, t_2, \ldots $ the following identities  are consequences of~(\ref{cma}):  
\begin{equation}
\label{Ja}
[x, y, t_1, \ldots, t_{2k}, z] +  [y, z, t_1, \ldots, t_{2k}, x] + [z, x, t_1, \ldots, t_{2k}, y] ,
\end{equation}
\begin{equation}
\label{C}
[x, y, t_1, \ldots, t_{k}, [u,v]] + (-1)^{m} [u, v, t_1, \ldots, t_{k}, [x,y]].
\end{equation}
\end{corollary}

\noindent {\bf Proof}. The identity~(\ref{Ja}) is an evident consequence of the Jacobi identity and Lemma~\ref{Ins}. 
The identity~(\ref{C}) for an even integer $k$  is obtained from~(\ref{Ja}) with the substitution $z\mapsto [u,v]$,
and for odd $k$ --- with the substitutions $z \mapsto [u,v]$, $y \mapsto t_k$, and $x \mapsto [x,y]$. \qed

\medskip

\begin{lemm}
 \label{cnz} 
The algebra  $M_{1,1}(E^1)$ does not satisfy the identity 
 $[x, y, t^{(m)}, [u,v]] = 0$ for any integer $m \ge 0$. 
\end{lemm}
\noindent  {\bf Proof.} The result is true because of the substitution
$x = u = t = \left(
\begin{array}{cc}
1 & 0 \\
0 & 0
\end{array}\right)$, $y = \left(
\begin{array}{cc}
0 & e_1 \\
0 & 0
\end{array}\right)$, and $v = \left(
\begin{array}{cc}
0 & 0 \\
e_2 & 0
\end{array}\right)$. (Here $e_1$ and $e_2$ are generators of $E$.) \qed

\begin{lemm}
\label{subs}
Let 
%$\mathbb F$ be a field and 
$f=0$ be a multihomogeneous identity of $M_{1,1}(E^{1})$ of degree $m \ge 3$. 
If   for some $k \ge 1$
$$f(\bar x) = \sum_{i, j \ge 2} \alpha_{i j} [x_1,x_i,\ldots , x_j] +\sum_{i \ge 1, i\ne k} \beta_{i } [x_k, x_i,\ldots , x_1]$$ 
then 
\begin{itemize}
\item[(i)] $\alpha_{ij} = (-1)^m \alpha_{ji}$;
\item[(ii)] $ \beta_j = (-1)^{m+1} \sum_{i} \alpha_{ij}$ for every  $j$ such that $j \ge 2$ and $j\ne k$;
\item[(iii)] $ \sum_{i} \beta_i = (-1)^m \sum_{i} \alpha_{ik}$ if $k>1$.
%\item[(4)] $ sum_{i} \beta_i = (-1)^m \sum_{i} \alpha_{ik}$ if $k>0$;
%\item[(2)] $2\beta_0 =  0$;  if $m$ is odd;
%\item[(5)] $ 2sum_{i} \beta_i = 0$ if $k=0$ and $m$ is odd.
\end{itemize}
\end{lemm}       

\noindent {\bf Proof}.
Denote by $A$ the commutator $[x,y]$ and by  $B$ the commutator $[u,v]$. 
To prove~(i)  substitute  in the identity $f=0$ the sum $x_i+A$ instead of $x_i$ and $x_j +B$ instead of 
 $x_j$, if $i\ne j$, and $x_i+A+B$ 
instead of $x_i$,  if $i=j$. Then one expands and denotes  by $g$ the sum of all summands which are not 
consequences of~$cm$ and contain 
both $A$ and $B$. It is clear that every summand from $g$ contains both $A$ and $B$ exactly once and 
the equality $g = 0$ is also an identity of $M_{1,1}(E^1)$ for any field $\mathbb F$. 
We have $g = \alpha_{i j} [x_1,A,\ldots , B] + \alpha_{j i} [x_1, B,\ldots, A]$. Let us use~(\ref{C}) to obtain
$$(\alpha_{i j} + (-1)^{m+1}\alpha_{j i}) [x_1,A,\ldots , B] = 0.$$
The latest equation and Lemma~\ref{cnz} provide~(i). 

For the case~(ii) we repeat the arguments above with the substitutions  $x_1 \mapsto x_1+A$ and 
$x_j \mapsto x_j + B$. To prove~(iii) we substitute $x_1+A$ instead $x_1$ and $x_k + B$
instead $x_k$ and also repeat the arguments of the previous  paragraph. \qed

%\medskip

%Now we are ready to prove Propositions~\ref{ml} and~\ref{mh1}.   

 \subsubsection*{Proof of Proposition~\ref{ml}}
%{\bf Proof of Proposition~\ref{ml}}. 

%The algebra $M_{1,1}(E)$ is a subalgebra of the algebra $M_{1,1}(E^1)$.  
%Notice that both $M_{1,1}(E^1)$  and $M_{1,1}(E)$. . We prove the theorem only for  $M_{1,1}(E)$. 
Let $f(x_1, x_2, \ldots, x_n)=0$ be a multilinear identity of $M_{1,1}(E)$. 
Then modulo~$cm$ and the Jacobi identity the polynomial $f$ can be written in the following form:
$$f = \sum_{i \ne j, \;  i, j \ge 2} \alpha_{ij} [x_1, x_i, \ldots, x_j].$$

Firstly, suppose that $n$ is odd. Then $n\ge 3$. By Lemma~\ref{subs}~(i) we have $\alpha_{ji} = - \alpha_{ij}$. 
Using~(\ref{Ja}) one can obtains 
\begin{align*}
f = & \sum_{2 \le i < j \le n} \alpha_{ij} ([x_1, x_i, \ldots, x_j] - [x_1, x_j, \ldots, x_i]) = 
 & \sum_{2 \le i < j \le n} \alpha_{ij} [x_j, x_i, \ldots, x_1]. 
\end{align*}
The Jacobi identity provides that $f= \sum_{i=3}^{n} \beta_{i} [x_2, x_i, \ldots, x_1]$. Hence, by Lemma~\ref{subs}~(i)
we have  $\beta_i = 0$ for any $i \ge 3$. This means that $f=0$ is a consequence of~$cm$. 

Now we assume that $m \ge 4$ is even. 
By~Lemma~\ref{subs}~(i)  
$$f =  \sum_{2 \le i < j \le n} \alpha_{ij} ([x_1, x_i, \ldots, x_j] + [x_1, x_j, \ldots, x_i]).$$
Notice that  if $i>2$ and $j>2$ we have 
\begin{align*}
& [\underbrace{x_1, x_i, \ldots ,x_2}_{m-1}, x_j] + [x_1, x_j, \ldots ,x_2, x_i] \stackrel{(\ref{Ja})}{=} \\
&-[\underbrace{x_2, x_1, \ldots ,x_i}_{m-1}, x_j] - [\underbrace{x_i, x_2, \ldots ,x_1}_{m-1}, x_j] + 
[x_1, x_j, \ldots ,x_2, x_i] = \\
& [x_1, x_2, \ldots ,x_i, x_j] + [x_2, x_i, \ldots , [x_1, x_j]] + [x_2, x_i, \ldots , x_j, x_1] + \\
&[x_1, x_j, \ldots, [x_2, x_i]] + [x_1, x_j, \ldots, x_i, x_2]  \stackrel{(\ref{C})}{=}  \\
 & [x_1, x_2, \ldots ,x_i, x_j] + [x_2, x_i, \ldots , x_j, x_1] + [x_1, x_j, \ldots, x_i, x_2]. 
\end{align*}
Hence, $f$ can be written in the following form
$$f = \sum_{i=3}^{n} \delta_{i}[x_1, x_2, \ldots, x_i] + \sum_{i=3}^{n} \gamma_{i} [x_1, x_i, \ldots, x_2] + 
\sum_{i=3}^{n} \beta_{i}  [x_2, x_i, \ldots, x_1].$$
By Lemma~\ref{subs}~(ii)  
we have $\gamma_i = \delta_i$ and $\beta_i = - \delta_i$ for every $i\ge 3$ and by Lemma~\ref{subs}~(iii) 
$\sum_{i} \beta_i = \sum_{i} \gamma_{i}$. Hence,  $2\sum_{i} \beta_i = 0$ and $\sum_{i} \beta_i = 0$.   
 It is clear that  $\beta_3 =  - \sum_{i=4}^{n} \beta_{i}$.  Therefore, we can rewrite $f$ 
in the following form:                                                                             
\begin{align*}
f = &\sum_{i=4}^{n} \beta_{i} ([x_1, x_2, \ldots, x_i] + [x_1, x_i, \ldots, x_2] - [x_2, x_i, \ldots, x_1]) - \\
( &\sum_{i=4}^{n} \beta_{i}) ([x_1, x_2, \ldots, x_3] + [x_1, x_3, \ldots, x_2] - [x_2, x_3, \ldots, x_1]) = \\
& \sum_{i=4}^{n} \beta_{i} ([x_1, x_2, \ldots, [x_3, x_i]] + [x_1, [x_i, x_3], \ldots, x_2] - [x_2, [x_i,x_3], \ldots, x_1]). 
\end{align*}
It remains to notice that for each $i$ the sum of these three elements in the brackets is a consequence of 
the identity~(\ref{Ja}).

 \subsubsection*{Proof of Proposition~\ref{mh1}}
%{\bf Proof of Proposition~\ref{mh1}}. 
Let $f$ be a multihomogeneous polynomial of degree $m$. One can assume that $f$ is not 
multilinear due to Theorem~\ref{ml}. Let $x_1$ be a variable from $f$ of degree greater than 1. 
Using the Jacobi identity one can write $f$ in the following form
$$f = \sum_{i \ne j , \; i, j \ge 2} \alpha_{ij} [x_1, x_i, \ldots, x_j] + \sum_{i \ge 2} \beta_{i} [x_1, x_i, \ldots, x_1].$$ 
If $m$ is odd then $\alpha_{ii}=0$ because of~(i) in Lemma~\ref{subs}. If $m$ is even then 
$[x_1, x_i, \ldots, [x_1,x_i]] = 0$, provided by~(\ref{C}). This means that one can replace  
$[x_1, x_i, \ldots, x_i]$ by $[x_1, x_i, \ldots, x_1]$. Anyway one can assume that  $\alpha_{ii}=0$.
 
By~(\ref{C}) the equality $[x_1, x_i, \ldots, [x_1,x_j]] = (-1)^{m+1}[x_1, x_j, \ldots, [x_1,x_i]]$ 
is a consequence of~(\ref{cma}). Hence, 
\begin{align*}
[x_1, x_i, \ldots, x_1,x_j] & = [x_1, x_i, \ldots, x_j,x_1] + \\
& (-1)^{m+1}[x_1, x_j, \ldots, x_1,x_i] + (-1)^{m+1}[x_1, x_j, \ldots, x_1,x_i].       
\end{align*}
Therefore, we can assume that $\alpha_{ij} = 0$ if $i \ge j$. Then by Lemma~\ref{subs}~(i) all $\alpha_{ij}$ 
equal zero. Now it remains to see that $\beta_i = 0$ for all $i$ because of Lemma~\ref{subs}~(ii).
This means that $f$ is a consequence of~$cm$.

%\subsubsection*{Proof of Theorems~~\ref{tinf} and\ref{t0}}
\medskip

\iffalse

Propositions imply Theorems~\ref{tinf} and~\ref{t0}. Indeed,   
%it is well-known that the identities of algebras over an infinite 
%field are generated by multihomogeneous identities. Thus, 
Theorem~\ref{tinf} is an evident consequence of 
Proposition~\ref{mh1} because of Remark (see page~\label{homog}). 

In the case of zero characteristic field identities of each algebra follows from multilinear ones. 
Hence, to show the truth of Theorem~\ref{t0} is sufficient to notice that the algebra 
$M_{1,1}(E^1)$ satisfies $cm=0$ and 
contains $M_{1,1}(E)$ as a subalgebra. This fact and Proposition~\ref{ml}  mean that 
$M_{1,1}(E^1)$ and $M_{1,1}(E)$ satisfy 
the same multilinear identities. Theorem~\ref{t0} is proved.    

\fi

 \subsection*{Proof of Theorem~\ref{Ep}}

The following result shows that for $M_{1,1}(E)$ the situation with  
identities over a finite field is not too different from another one, when the base field is infinite.

\begin{lemm}
\label{ff}
Let $\mathbb F$ be a finite field. Then the multihomogeneous components 
of every identity of $M_{1,1}(E)$  are also identities of the algebra.    
\end{lemm}

\noindent {\bf Proof.} Assume the contrary. Then there exists $f(x, \ldots )$, an identity   
of $M_{1,1}(E)$, equal to the sum of non-identities 
$g(x,\ldots )$ and $h(x,\ldots )$, where $g$ is homogeneous in $x$ of degree $n$ and 
%a variable $x$ includes in 
%each monomial of $g$ exactly $n$ times but in each monomial of $h$ less than $n$ times. 
each monomial of $h$ has degree in $x$ less than $n$.  
One can suppose that $n$ is the minimal integer greater than one with the condition. By Remark (see p.~\pageref{homog}), 
we conclude that $n \ge p$, 
if $p$ is a characteristic of $\mathbb F$. 
%(Otherwise, we can apply the standard method using Vandermonde determinant 
%to show that $g$ and $h$ are also identities.) 
Let  $\tilde f$, $\tilde g$, $\tilde h$ are results of the linearization of the polynomials, 
that is we have for $f$  
$$\tilde f(x,y, \ldots) = f(x+y, \ldots)- f(x \ldots) - f(y, \ldots),$$
and the similar equalities for $\tilde g$ and $\tilde h$. 
As far as $n$ is supposed to be minimal we conclude that $\tilde g$ is an identity of $M_{1,1}(E)$. 
Indeed, if $n=2$, then $\tilde g$ is an identity because then $deg_x (h) = 1$ 
%includes in monomials of $h$ only once 
and 
$\tilde h$ becomes zero. If $n>2$, then for each $i=1,\ldots, n-1$ we denote by 
$\tilde g_i$, the homogeneous in $x$ component of degree $i$ in $x$ from $\tilde g$. 
% sum all monomials of degree $i$ in $x$ from $\tilde g$. 
The polynomial $\tilde g_{n-1}$ is an identity  
of our algebra because otherwise one can consider it as a ``new'' polynomial $g$ having a degree in $x$ less than $n$. 
Clearly, $\tilde g_1$ is  also an identity because it is obtained from $\tilde g_{n-1}$ by the 
transposition of $x$ and $y$.
Repeating the argument we show that all  $\tilde g_i$ are identities of the algebra as well as  $\tilde g$. 

Anyway, $\tilde g$ is an identity of $M_{1,1}(E)$. This means that $g$ can be considered as a map 
acting on our algebra and linear in $x$. It is easy to see that each  
element from $M_{1,1}(E)$ is a sum of such elements $a_i$ that 
an arbitrary monomial $w$ with $deg_{t}(w) \ge p$ becomes zero after the substitution $t \mapsto a_i$, whatever the choice
of $i$. 
%an arbitrary word containing $p$ copies of $a_i$ for some $i$  equals zero. 
As far as $n>p$ we have that  $g$ becomes 
zero for every substitution of elements from our algebra. A contradiction.   \qed

\medskip         

Now we need some technical results and notations. Denote

\begin{gather}
\label{cp}
[x,z,u^{(p)},v]; \\
\label{pp}
[x,y,x^{(p-1)},y^{(p-1)},v].
\end{gather}

%$\mathcal{V} = \var\{ [x,y,[z,u],v], [x,z,u^{(p)},v], [x,y,x^{(p-1)},y^{(p-1)},v] \}$

\begin{lemm}
\label{conseq}
If $\;\mathbb F$ is a field of characteristic $p>2$, then for each positive integer $k$ the following identities 
are consequences of~(\ref{cma}),~(\ref{cp}), and~(\ref{pp}):
\begin{gather}
\label{con1}
[y,z_1, \ldots, z_{2k}, y^{(p)}]; \\
\label{con2}
[x,y^{(p)},z_1, \ldots, z_{2k-1}, t] + [x, t ,z_1, \ldots, z_{2k+1}, y^{(p)}] - [x, y, t, z_1, \ldots, z_{2k+1},y^{(p-1)}];\\
\label{con3}
[x,y^{(p)},z_1\ldots ,z_{2k-1}, x^{(p-1)}] - [y,x^{(p)},z_1\ldots ,z_{2k-1}, y^{(p-1)}].
\end{gather}
\end{lemm}

\noindent  {\bf Proof}. 
We have  
\begin{align*}
& 0 \stackrel{(\ref{cp})}{=} [z_2,z_1, y^{(p)}, \ldots, x] \stackrel{(\mbox{\tiny Jacobi identity})}{=}  - [z_1,y, z_2, y^{(p-1)}, \ldots,x]  -  
[ [y, z_2], z_1, y^{(p-1)}, \ldots,x]  \stackrel{(\ref{Ja})}{=} \\
 - &[z_1,y, z_2, y^{(p-1)}, \ldots,x]  +  [ x, [y, z_2], y^{(p-1)}, \ldots, z_1] + [ z_1, x,  y^{(p-1)}, \ldots, [y, z_2]] = \\
 - &[z_1,y, z_2, y^{(p-1)}, \ldots,x]  -  [ y, z_2, x, y^{(p-1)}, \ldots, z_1] +  \\ 
& [ z_1, x,  y^{(p-1)}, \ldots, y, z_2] - [ z_1, x,  y^{(p-1)}, \ldots, z_2, y] \stackrel{(\ref{cma}), (\ref{cp})}{=} \\ 
 - &[z_1,y, z_2, y^{(p-1)}, \ldots,x]  -  [ y, z_2, x, y^{(p-1)}, \ldots, z_1] - [ z_1, x,  y^{(p-1)}, \ldots, z_2, y]. 
\end{align*}
%The third summand equals zero modulo~(\ref{cma}) and~(\ref{cp}). The second summand becomes the same 
The second summand becomes zero modulo~(\ref{cma}) and~(\ref{cp}) 
after the substitution $x \mapsto y$. Hence, after the substitution we have
 modulo~(\ref{cma}) and~(\ref{cp})   $$0 = 2[y, z_1, z_2,  \ldots, y^{(p)}].$$ This  proves 
that~(\ref{con1}) follows from~(\ref{cma}) and~(\ref{cp}).

 In the case of infinite field $\mathbb F$ the identity~(\ref{con2}) is a partial linearization of~(\ref{con1}). But we obtain 
the identity for an arbitrary field. At first notice that due to the Jacobi identity the following equality holds modulo~(\ref{cp}):
\begin{equation}
\label{tr}
[x,y^{(p)},z,v] = [z,y^{(p)},x,v].
\end{equation}
Moreover, we have
\begin{align*}
& 0 \stackrel{(\ref{C})}{=} [x,y, y^{(p-2)},z_1, \ldots, z_{2k-1}, [y, t]] + [y, t, y^{(p-2)},z_1, \ldots, z_{2k-1}, [x,y]] = \\  
%& [x,y, y^{(p-2)},z_1, \ldots, z_{2k-1}, y, t] - [x,y, y^{(p-2)},z_1, \ldots, z_{2k-1}, t, y] + [y, t, y^{(p-2)},z_1, \ldots, z_{2k-1}, x,y] - [y, t, y^{(p-2)},z_1, \ldots, z_{2k-1}, y, x] = \\ 
& [x,y, \ldots,  y, t] - [x,y, \ldots,  t, y] + [y, t,  \ldots,  x,y] - [y, t,  \ldots,  y, x] \stackrel{(\mbox{\tiny Jacobi identity}), (\ref{tr})}{=}\\ 
& [z_1,y^{(p)}, x \ldots, t] + [t, x, \ldots,  y^{(p)}] + [y, t, \ldots,  x,y] + [y, t, \ldots,  x,y] + [z_1, y^{(p)}, t,  \ldots, x] = \\ 
& 2[z_1,y^{(p)}, x \ldots, t] + [t, x, \ldots,  y^{(p)}] + 2[y, t, \ldots,  x,y] +  [z_1, y^{(p)},  \ldots, [t,x]] \stackrel{(\ref{cp})}{=}\\ 
& 2[z_1,y^{(p)}, x \ldots, t] + [t, x, \ldots,  [z_1, y^{(p)}]] + [z_1, y^{(p)},  \ldots, [t,x]] + 2[y, t, \ldots,  x,y]   \stackrel{(\ref{C})}{=}\\ 
& 2([z_1,y^{(p)}, x \ldots, t] + [y, t, \ldots,  x,y]) \stackrel{\mbox{\tiny (Jacobi identity)}}{=}\\ 
& 2([x,y^{(p)}, z_1, \ldots, t]  - [x, y, \ldots,  t,y] + [t, x,\ldots , y^{(p)}]). 
\end{align*}

Let us prove~(\ref{con3}).
\begin{align*}
0 & \stackrel{(\ref{pp})}{=}  [x,y, x^{(p-1)}, y^{(p-1)}, z_1,\ldots , z_{2k-1}] \stackrel{(\ref{Ja})}{=} \\
&  -  [ y, z_{2k-1}, x^{(p-1)}, y^{(p-1)}, z_1,  \ldots , z_{2k-3}, x] - [z_{2k-1}, x^{(p)}, y^{(p)}, z_1, \ldots , z_{2k-3}, y] 
\stackrel{(\ref{tr}), (\ref{cma})}{=} \\ 
& [x, y^{(p)}, z_1,\ldots , z_{2k-1}, x^{(p-1)}] - [y, x^{(p)}, z_1, \ldots ,z_{2k-1}, y^{(p-1)}]. 
\end{align*}

\begin{lemm}
\label{nonc}
Let $\mathbb F$ be a field of characteristic $p>2$. The algebra $M_{1,1}(E)$ does not satisfy the  
following identities for any integer $k\ge 0$:
$$[y,x, z_1 \ldots, z_{2k}, y^{(p)}] = 0;$$
$$[x,y^{(p-1)},x^{(p-1)}, z_1\ldots ,z_{2k}, y] = 0.$$
\end{lemm}
\noindent  {\bf Proof.}  Let us put 
$z_i =  \left(
\begin{array}{cc}
a_i & 0 \\
0 & 0
\end{array}\right)$, 
$x = \left(
\begin{array}{cc}
b_0 & b_1 \\
0 & 0
\end{array}\right)$, and $y = \left(
\begin{array}{cc}
c_0 & c_1 \\
c_2 & 0
\end{array}\right)$. Here $a_i$, $b_0$, $c_0$ are from $E_0$, and $b_1, c_1, c_2$ are from $E_1$.
Then $$[y,x, z_1 \ldots, z_{2k}, y^{(p)}] =   \left( \begin{array}{cc}
w & 0 \\
0 & w
\end{array}\right),$$ 
where $w=2a_1\cdots a_{2k}b_0c_0^{p-1}c_1c_2 $. Clearly, $w$ is not an identity of $E$. 
After the substitutions 
the element $$[x,y^{(p-1)},x^{(p-1)}, z_1\ldots ,z_{2k}, y] $$ turns into the 
diagonal matrix with the element 
$a_1\cdots a_{2k}b_0^{p-1}c_0^{p-1}b_1c_2$ at the diagonal. \qed

\medskip

\noindent {\bf Proof of Theorem~\ref{Ep}}.
It is easy to show that~(\ref{cp}) and~(\ref{pp}) are identities of $M_{1,1}(E)$.

Now we want to make sure that every identity of $M_{1,1}(E)$ follows from~(\ref{cma}),~(\ref{cp}) 
and~(\ref{pp}). Due to Lemma~\ref{ff} we can consider only multihomogeneous polynomial identities of our algebra. 
Let $f$ be a such  identity  of  degree $m$.	
If some variable in $f$ has a degree less than $p$, then one can linearize $f$ in the variable to obtain 
the polynomial which $f$ follows from, and which is a sum of multihomogeneous identities. 
The variable includes in each such summand 
in degrees less than in $f$. Thus,  we can assume that every variable in $f$ has either a degree equal to 1 or a 
degree greater than or equal to $p$.    
By Proposition~\ref{ml} one can assume that $f$ is not multilinear. 
Let us denote by $y$ the variable from $f$ of the greatest degree and by $x$ the greatest degree variable among 
the variables without $y$.  It is clear that the degree of $y$ is 
greater than or equal to $p$. We assume that $f$ is not a consequence of the identities~(\ref{cma}) 
and~(\ref{cp}). This means in particular that the degree of $y$ is less than $p+2$. 
 Using the Jacobi identity and~(\ref{cma}) one can 
write the polynomial $f$ in the following form: 
$$f= \sum_{i}\alpha_{i}[x,y,y, \ldots, x_i] + \sum_{i}\beta_{i}[x, x_i,\ldots,y,y] + \gamma [x,y, \ldots, y],$$       
where the variables $x_i$ don't coincide with $y$ but can coincide with $x$. 

{\bf Case 1 ($deg_y f =p + 1$).} It is evident that except for the last
summand of $f$, all of them are consequences of~(\ref{cma}) and~(\ref{cp}). 
Hence, one can assume that  $f = \gamma [x,y, \ldots, y]$.
If $m$ is even then $f=0$ is a consequence 
of~(\ref{con1}). If $m$ is odd and the degree of $x$ is equal to 1, then
$  \gamma = 0$ by Lemma~\ref{nonc}. 
If $m$ is odd and the degree of $x$ is  greater than 1, then $f=0$ is an identity of $\mathcal{V}$.  
 
{\bf Case 2 ($deg_y f =p $, $deg_x f =1$).} 
By Lemma~\ref{nonc} $f$ depends on greater than 2 variables.

Let us substitute $x_i \mapsto A$, where $A=[u,v]$. Then modulo identities~(\ref{cma}) and~(\ref{cp}) we have 
\begin{align*}
0 = & \alpha_{i}[x,y,y, \ldots, A] + \beta_{i}[x, A,\ldots,y,y] \stackrel{(\ref{C})}{=} \\
 & (-1)^m\alpha_{i}[A,y, \ldots, [x,y]] - \beta_{i}[A,x, \ldots,y,y] \stackrel{(\ref{cma})}{=} \\
 & (-1)^m\alpha_{i}[A,x, \ldots, y,y] - (-1)^m\alpha_{i}[A,y^{(p)}, \ldots, x] - \beta_{i}[A,x, \ldots,y,y] \stackrel{(\ref{cp})}{=}\\ 
& ((-1)^m\alpha_{i} - \beta_{i}) [A,x, \ldots,y,y] 
\end{align*}
Hence, we have $(-1)^m\alpha_{i} = \beta_{i}$ provided by Lemma~\ref{nonc} for every $i$.

If $m$ is odd then 
\begin{align*}
f = & \sum_i \alpha_{i}([x,y,y, \ldots, x_i] - [x, x_i,\ldots,y,y])  + \gamma [x,y, \ldots, y]\stackrel{(\ref{Ja}),(\ref{cma}) }{=} \\
 &  \sum_i \alpha_{i}[x_i,y,y, \ldots, x] +  \gamma [x,y, \ldots, y]  \stackrel{\mbox{(\tiny Jacobi identity)}, (\ref{cp})}{=} \\
 & \delta [x_1,y,y, \ldots, x] +  \gamma [x,y, \ldots, y].
\end{align*}
It remains to substitute $x\mapsto y$ to obtain $\delta = 0$ provided by Lemma~\ref{nonc}. 
Then this lemma guarantees that
$\gamma = 0$.     

If $m$ is even then 
\begin{align*}
f = & \sum_i \alpha_{i}([x,y,y, \ldots, x_i] + [x, x_i,\ldots,y,y])  + \gamma [x,y, \ldots, y] \stackrel{(\ref{con2}),(\ref{cma}) }{=} \\
 &  \sum_i \alpha_{i}[x,y, \ldots, y]  +  \gamma [x,y, \ldots, y] = \delta [x,y, \ldots, y].
\end{align*}
In this case we substitute $x\mapsto [x,y]$ to use Lemma~\ref{nonc} and obtain $\delta = 0$.

{\bf Case 3 ($deg_y f = p $, $deg_x f = p $).} 

First we notice that  if  $x_i \ne x$ then $[x, y,\ldots, x_i]=0$ is a consequence of~(\ref{pp}). Furthermore, 
by~(\ref{C}) we have  
\begin{align*}
& [x, x_i,\ldots,x,y] = (-1)^{m} [x, y,\ldots,[x,x_i]] + [x, x_i,\ldots, y, x]  =  \\
& (-1)^{m} [x, y,\ldots, x_i] - (-1)^{m} [x, y,\ldots, x_i, x] + [x, x_i,\ldots, y, x] 
\stackrel{(\ref{cp}), (\ref{pp}), (\ref{cma}) }{=} \\
 &  - (-1)^{m} [x, y,\ldots, x_i, x].
\end{align*}
These two facts above mean that 
$$f = \phi [x,y, \ldots, x]  + \gamma [x,y, \ldots, y].$$
If $m$ is even then one can use~(\ref{con2}) with $t = x$ to obtain $$[x,y, \ldots, x] =  [x,y, \ldots, y].$$
The same result implies from~(\ref{con3}) for odd $m$. Anyway, 
 $f = (\phi + \gamma) [x,y, \ldots, y]$. By Lemma~\ref{nonc} we have $\phi + \gamma = 0$. Theorem~\ref{Ep} is proved.


\begin{thebibliography}{99}

\bibitem{AK06}  S.M. Alves, P. Koshlukov, Polynomial identities of algebras in positive characteristic, 
J. Algebra {\bf 305} (2) (2006) 1149–-1165.

\bibitem{AFK04} S. Azevedo, M. Fidelis, P. Koshlukov, Tensor product theorems in positive characteristic, J. Algebra, 
{\bf 276} (2), (2004) 836–-845.

\bibitem{AFK05} S.S. Azevedo, M. Fidelis, P. Koshlukov, Graded identities and PI equivalence of algebras in positive characteristic, Comm.
Algebra {\bf 33} (4) (2005) 1011–-1022.

\bibitem{Kemer91} A. Kemer, Ideals of identities of associative algebras. Translations Math.
Monographs 87. Providence, RI: Am. Math. Soc. (1991).

\bibitem{KM}  P. Koshlukov, T.C. de Mello, On the polynomial identities of the algebra $M_{1,1}(E)$,  
Linear Algebra Appl. {\bf 438} (11) (2013), 4469--4482. 

\bibitem{Kras2008} A. N. Krasilnikov, The identities of a Lie algebra viewed as a Lie ring, Quart. J.
Math., {\bf 60} (1) (2009), 57–-61.

\bibitem{Kuzmin} Yu. V. Kuz'min, Free center-by-metabelian groups, Lie algebras, and D-groups, Math. USSR-Izv., {\bf 11}, 
1 (1977), 1--30.

\bibitem{Mishch} S. P. Mishchenko, Varieties of centrally metabelian Lie algebras over a field of characteristic zero, 
Math. Notes, 305 (1981), 822-–826. 

\bibitem{Popov} A. Popov, Identities of the tensor square of a Grassmann algebra, Algebra Log. {\bf 21} (4) 
(1982), 296-–316. 
%Algebra i Logika, 21 (4) (1982), pp. 442-471

\bibitem{Samoilov} L. M. Samoilov, On multilinear components of prime subvarieties in the variety $Var(M_{1,1})$, 
Math Notes, {\bf 87}, (2010), 890--902. 
%Matematicheskie Zametki, {\bf 87} (6) (2010),  919--933.

\bibitem{Vincenzo} O. M. Di Vincenzo,  On the graded identities of $M_{1,1}(E)$,  Israel J. Math. {\bf 80} (3) (1992),  323--335. 

\bibitem{Voliche} I. B. Volichenko, Varieties of centre-by-metabelian Lie algebras. 
Prepr. 1696, Inst. Mat. AN BSSR (1980b) (Russian).

\end{thebibliography}
\end{document}